\documentclass[10pt]{article}
\usepackage{amsmath, amssymb, amsthm, epsfig, color}

\newtheorem{theorem}{Theorem}

\newtheorem{lemma}{Lemma}

\theoremstyle{definition}
\newtheorem{remark}{Remark}
\newtheorem{example}{Example}
\newtheorem{definition}{Definition}

\title{Path counting and random matrix theory}
\author{Ioana Dumitriu and Etienne Rassart\footnote{Supported by FCAR (Qu\'ebec)}\\[2mm] \emph{Department of Mathematics, Massachusetts Institute of Technology}}

\begin{document}
\maketitle

\begin{abstract} 
We establish three identities involving Dyck paths and alternating Motzkin paths, whose proofs are based on variants of the same bijection. We interpret these identities in terms of closed random walks on the halfline. We explain how these identities arise from combinatorial interpretations of certain properties of the $\beta$-Hermite and $\beta$-Laguerre ensembles of random matrix theory. We conclude by presenting two other identities obtained in the same way, for which finding combinatorial proofs is an open problem.
\end{abstract}

\section{Overview}
In this paper we present five identities involving Dyck paths and alternating Motzkin paths. These identities appear as consequences of algebraic properties of certain matrix models in random matrix theory, as briefly described in Section \ref{defs}. Three of them describe statistics on Dyck and alternating Motzkin paths: the average norm of the rise-by-altitude and vertex-by-altitude vectors for Dyck paths, and the weighted average square norms of the rise-by-altitude and level-by-altitude vectors for alternating Motzkin paths. We describe these quantities in detail in Section \ref{defs}, and provide combinatorial proofs for the identities in Section \ref{the3}.

In terms of closed random walks on the halfline, these identities give exact formulas for the total square-average time spent at a node, as well as the total square-average number of advances to a higher labelled node.

For the other two identities we have not been able to find simple interpretations or combinatorial proofs that would complement the algebraic ones; this is a challenge that we propose to the reader in Section \ref{the2}.

\section{Definitions, main results, and interpretations} \label{defs}

The Catalan numbers $C_k$ count dozens of combinatorial structures, from binary trees and triangulations of polygons to Dyck paths \cite[Exercise 6.19, pages 219-229]{EC2}. Similar, but less known, are the Narayana numbers $N_{k,r}$ \cite[Exercise 6.36, page 237]{EC2}; since they sum up to $C_k$, they partition combinatorial structures enumerated by Catalan numbers according to a certain statistic. In particular, they count alternating Motzkin paths (see Section \ref{the3}). 

The relationship between Catalan numbers and random matrix theory appeared first in Wigner's 1955 paper \cite{wigner55a}. In computing asymptotics of traces of powers of certain random (symmetric, hermitian) matrices, Wigner obtained (not explicitly by name) the Catalan numbers, which he recognized as the moments of the semi-circle. Later, Mar{\v{c}}enko and Pastur, in their 1967 paper \cite{marcenko67a} found a similar connection between Narayana numbers and Wishart (or Laguerre) matrix models (more explicitly, they computed the generating function for the Narayana polynomial). 
Both connections have to do with computing average traces of powers of random matrices, i.e. the moments of the eigenvalue distribution. 

Suppose $A$ is an $n \times n$ symmetric random matrix, scaled so that as $n \rightarrow \infty$ the probability that its eigenvalues lie outside of a compact set goes to $0$. Denoting by  
\[
m_k = \lim_{n \rightarrow \infty} E\left[ \frac{1}{n} \mbox{tr}(A^k) \right]~,
\] 
one can ask the question of computing $m_k$ for certain types of random symmetric matrix models. In some cases, $m_k$ might not even exist, but in the cases of the Gaussian and Wishart matrix models, it does. For the Gaussian model, 
\[
m_k = \left \{ \begin{array}{ll} 0, & \mbox{if $k$ is odd}, \\
				C_{k/2}, & \mbox{if $k$ is even}. \end{array} \right. ~,
\]
while for the Wishart model $W = GG^{T}$, where $G$ is a rectangular $m \times n$ matrix of i.i.d. Gaussians,
\[
m_k = N_k(\gamma)~,
\]
provided that $m/n \rightarrow \gamma$.

In both cases, one way of computing the zeroth-order (i.e. asymptotically relevant) term in $E\left[ \frac{1}{n} \mbox{tr}(A^k) \right]$ is by writing
\begin{eqnarray} \label{full_sum}
\mbox{tr}(A^k) = \sum_{i = 1}^{n} \sum_{1 \leq i_1, \ldots, i_{k-1} \leq n} a_{i i_1} a_{i_1 i_2} \ldots a_{i_{k-2}i_{k-1}} a_{i_{k-1} i}~,
\end{eqnarray}
then identifying the asymptotically relevant terms, weighing their contributions, and ignoring the rest. For example, if $k$ is even, in the case of the Gaussian models (which have i.i.d. Gaussians on the off-diagonal, and i.i.d. Gaussians on the diagonal), the only terms $a_{i i_1} \ldots a_{i_{k-1} i}$ which are asymptotically relevant come from sequences $i_0 = i, i_1, \ldots, i_{k} = i$ such that each pair $i_j, i_{j+1}$ appears exactly once in this order, and exactly once reversed. The connection with the Catalan numbers becomes apparent, as the problem reduces thus from counting closed random walks of length $k$ on the complete graph (with loops) of size $n$, to counting plane trees with $k/2$ vertices.

The above assumes full matrix models $A$; using the (equivalent) tridiagonal matrix models $T$ associated with a larger class of Gaussian and Wishart ensembles described in \cite{dumitriu02}, we can replace the problem of counting closed random walks on the complete graph to counting closed random walks on a line.

Using the tridiagonal model simplifies (\ref{full_sum}) to
\begin{eqnarray} \label{expansion}
\mbox{tr}(T^k) = \sum_{i = 1}^{n} \sum_{1 \leq i_1, \ldots, i_{k-1} \leq n} t_{i i_1} t_{i_1 i_2} \ldots t_{i_{k-2}i_{k-1}} t_{i_{k-1} i}~,
\end{eqnarray}
where $t_{i_j i_{j+1}}$ is non-zero iff $|i_j - i_{j+1}| \in \{0, \pm 1\}$. These correspond to closed walks on the line with loops. 

For the Gaussian models, when $k$ is even, the only asymptotically relevant terms can be shown to be given by closed walks which use no loops, which are in one-to-one correspondence with the Dyck paths of length $k/2$. For the Wishart models, these are closed walks on the line with loops that go right only on even time-steps, and left only on odd time-steps. In turn, these are in one-to-one correspondence with the alternating Motzkin paths.

The connection between Dyck paths, alternating Motzkin paths, and random matrix theory can be pushed further. In computing the variance of the traces of these powers for the Hermite and Laguerre ensembles, it can be shown algebraically \cite{dumitriu03a} that the zeroth and first-order terms in $n$ disappear. When one examines the expansion \eqref{expansion} applied to the tridiagonal models for Hermite and Laguerre ensembles, this translates into Theorems \ref{id1}, \ref{id2}, and \ref{id3}.

First, we recall the definitions of Catalan and Narayana numbers.

\begin{definition} The $k$th Catalan number $C_k$ is defined as
\[
C_k  = \frac{1}{k+1} {2k \choose k}~.
\]
\end{definition}

\begin{definition} The $(k,r)$ Narayana number is defined as 
\[
N_{k,r} = \frac{1}{r+1} {k \choose r} {k-1 \choose r}~.
\]
The associated Narayana polynomial (or generating function) is defined as
\[
N_k(\gamma) \equiv \sum_{r=0}^{k-1} \gamma^r N_{k,r} = \sum_{r=0}^{k-1} \gamma^{r} ~\frac{1}{r+1} {k \choose r} {k-1 \choose r}~.
\]
Note that $N_k(1)=C_k$.
\end{definition}

The Catalan numbers count many different combinatorial structures; in particular, they count Dyck paths.

\begin{figure}[ht]
\begin{center}
\epsfig{figure = 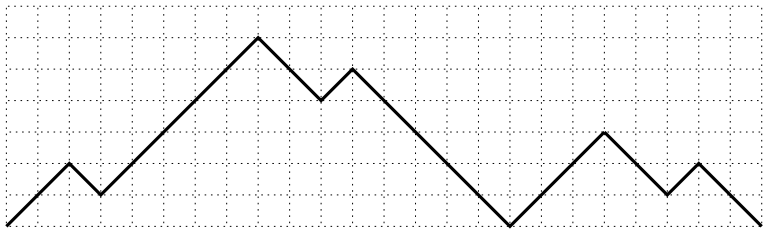, height = 3cm}
\end{center}
\caption{A Dyck path of length 24.} \label{DyckPath}
\end{figure}

\begin{definition}
A Dyck path of length $2k$ is a lattice path consisting of ``rise'' steps or ``rises'' ($\nearrow$) and ``fall'' steps or ``falls'' ($\searrow$), which starts at $(0,0)$ and ends at $(2k,0)$, and does not go below the $x$-axis (see Figure \ref{DyckPath}). We denote by $\mathcal{D}_k$ the set of Dyck paths of length $2k$.
\end{definition}

The Narayana numbers $N_{k,r}$ count alternating Motzkin paths of length $2k$ with $r$ rises; we recall the definition of Motzkin paths and define alternating Motzkin paths below.


\begin{definition} A Motzkin path of length $2k$ is a path consisting of ``rise'' steps or ``rises'' ($\nearrow$), ``fall'' steps or ``falls'' ($\searrow$), and ``level'' steps ($\rightarrow$), which starts at $(0,0)$, ends at $(2k,0)$, and does not go below the $x$-axis. 
\end{definition}

\begin{figure}[ht]
\begin{center}
\epsfig{figure = 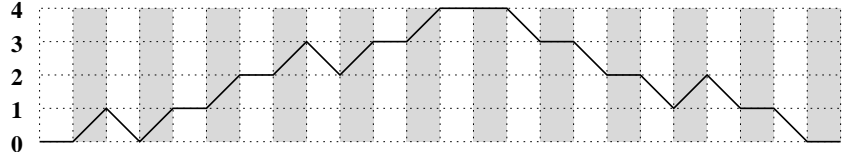, height= 2cm}
\end{center}
\caption{An alternating Motzkin path of length 24, with a total of 7 rises.} \label{AltDyckPath}
\end{figure}

\begin{definition} An alternating Motzkin path of length $2k$ is a Motzkin path in which rises are allowed only on even numbered steps, and falls are only allowed on odd numbered steps. See Figure \ref{AltDyckPath}. We denote by $\mathcal{AM}_k$ the set of alternating Motzkin paths of length $2k$.
\end{definition}

\begin{remark}
It follows from the definition that an alternating Motzkin path starts and ends with a level step.
\end{remark}

Next, we introduce three statistics on Dyck and alternating Motzkin paths.

\begin{definition} Let $p$ be a Dyck or alternating Motzkin path of length $2k$. We define the vectors $\vec{R} = \vec{R}(p) = (R_0, R_1, \ldots, R_{k-1})$ and $\vec{V} = \vec{V}(p) = (V_0, V_1, \ldots, V_{k})$ to be the rise-by-altitude and vertex-by-altitude vectors, i.e. $R_i$ is the number of rises from level $i$ to level $i+1$ in $p$, and $V_i$ is the number of vertices at level $i$ in $p$.
\end{definition}

For example, for the Dyck path of Figure \ref{DyckPath}, for which $k = 12$, 
\begin{eqnarray*}
\vec{R} & = & (2, 4, 2, 1, 2, 1, 0, 0, 0, 0, 0, 0)~, \\
\vec{V} & = & (3, 6, 6, 3, 3, 3, 1, 0, 0, 0, 0, 0, 0)~.
\end{eqnarray*}

Note that for a Dyck path of length $2k$, $\sum_{i=0}^{k-1} R_i = k$, while $\sum_{i=0}^{k} V_i = 2k+1$. For an alternating Motzkin path of length $2k$ with $r$ rises, $\sum_{i=0}^{k-1} R_i = r$, while $\sum_{i=0}^{k} V_i = 2k+1$.

\begin{definition} Let $p$ be an alternating Motzkin path of length $2k$. We define the vector $\vec{L} = \vec{L}(p) =  (L_0, L_1, \ldots, L_{k-1})$ be the even level-by-altitude vector, i.e. $L_i$ is the number of level steps at altitude $i$ in $p$ which are on even steps.
\end{definition}

\begin{remark} 
In the closed walk on a line interpretation, a rise from altitude $i$ to level $i+1$ corresponds to entering node $i+1$ from the left; a level step at altitutde $i$ corresponds to a loop from node $i$, and the number of vertices at altitude $i$ counts the number of time-steps when the walk is at node $i$.
\end{remark}

We are now able to state the three results, proved in Section \ref{the3}.

\begin{theorem} \label{id1} Let $F_{\mathcal{D}_k}$ be the uniform distribution on the set of Dyck paths of length $2k$. Then 
\[
\| E[\vec{R}]\|_2^2 ~~\equiv~~ \frac{1}{C_k^2} ~\sum_{p_1,p_2 \in \mathcal{D}_k} \sum_{i=0}^{k-1} R_i(p_1)R_i(p_2)~~ =~~\frac{C_{2k}}{C_k^2} - 1~,
\]
where $E$ denotes expectation with respect to $F_{\mathcal{D}_k}$.
\end{theorem}

\begin{remark}
In the closed random walk on the halfline interpretation, this identity gives a closed form for the total square-average number of advances to a higher labelled node.
\end{remark}

\begin{example}
\label{example_id1}
Here is an example for $k = 3$ of computing the average rise-by-altitude vector $\vec{R}$ and the average vertex-by-altitude vector $\vec{V}$ for Dyck paths of length 6. 

\begin{displaymath}
\begin{array}{c@{\hspace{8mm}}c@{\hspace{8mm}}c@{\hspace{8mm}}c@{\hspace{8mm}}c}
\epsfig{file=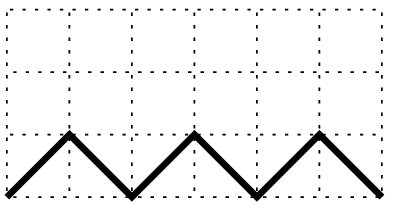, width=0.15\textwidth} & \epsfig{file=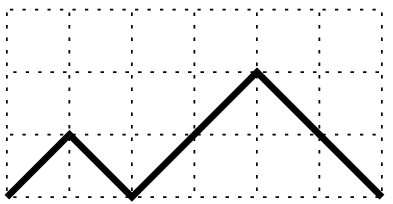, width=0.15\textwidth} & \epsfig{file=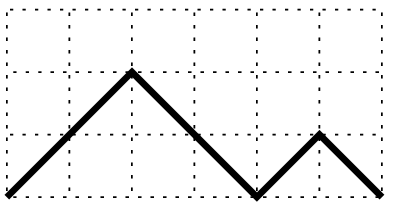, width=0.15\textwidth} & \epsfig{file=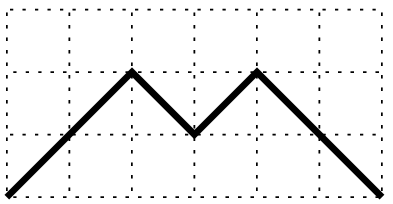, width=0.15\textwidth} & \epsfig{file=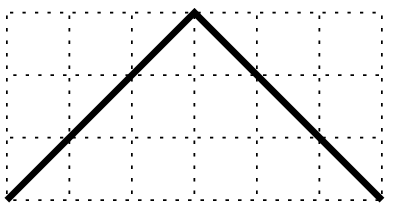, width=0.15\textwidth}\\
\vec{R}=(3,0,0) & \vec{R}=(2,1,0) & \vec{R}=(2,1,0) & \vec{R}=(1,2,0) & \vec{R}=(1,1,1)\\
\vec{V}=(4,3,0,0) & \vec{V}=(3,3,1,0) & \vec{V}=(4,3,0,0) & \vec{V}=(2,3,2,0) & \vec{V}=(2,2,2,1)
\end{array}
\end{displaymath}
\begin{displaymath}
E[\vec{R}] = \frac{1}{5}(9,5,1) \qquad \qquad E[\vec{V}] = \frac{1}{5}(14,14,6,1)
\end{displaymath}

Hence, for $k = 3$, 
\[
\| E[\vec{R}]\|_2^2 = \frac{81+25+1}{25} = \frac{107}{25} = \frac{C_6}{C_3^2} - 1~.
\]
\end{example}

\begin{theorem} \label{id2}
 Let $F_{\mathcal{D}_k}$ be the uniform distribution on the set of Dyck paths of length $2k$. Then 
\[
\| E[\vec{V}]\|_2^2 ~~\equiv ~~ \frac{1}{C_k^{2}} \sum_{p_1,p_2 \in \mathcal{D}_k} \sum_{i=0}^{k} V_i(p_1)V_i(p_2) ~~=~~ \frac{C_{2k+1}}{C_k^2}~,
\]
where $E$ denotes expectation with respect to $F_{\mathcal{D}_k}$.
\end{theorem}

\begin{remark}
In the closed random walk on the halfline setup, this gives a closed form for the total square-average time spent at a node.
\end{remark}

We use once again Figure \ref{example_id1}; 
\[
\| E[\vec{V}]\|_2^2 ~= \frac{196+196+36+1}{25} = \frac{429}{25} = \frac{C_7}{C_3^2}~.
\]

Finally, the third main result.

\begin{theorem} \label{id3}
Let $\gamma>0$, and let $F_{\mathcal{AM}_{k}}(\gamma)$ be the distribution on $\mathcal{AM}_{k}$ which associates to each alternating Motzkin path $p$ a probability proportional to $\gamma^r$, where $r$ is the number of rises in $p$. Then 
\begin{eqnarray*}
\| E[\vec{R}]\|^2_2~+~\gamma~\| E[\vec{L}]\|_2^2 & \equiv & \frac{1}{N_k(\gamma)^2} \sum_{p_1,p_2 \in \mathcal{AM}_k} \gamma^{r_1+r_2} \left(\sum_{i=0}^{k-1} R_i(p_1)R_i(p_2) + \gamma \sum_{i=0}^{k-1} L_i(p_1)L_i(p_2) \right) \\
 & = & \frac{N_{2k}(\gamma)}{N_k(\gamma)^2} - 1~,
\end{eqnarray*}
where $r_1$ and $r_2$ are the number of rises in $p_1$ and $p_2$, and $E$ denotes expectation with respect to $F_{\mathcal{AM}_{k}}(\gamma)$. 
\end{theorem}

\begin{remark}
In the closed random walk on the halfline setup, this gives a
relationship between the total square-average number of advances to a
higher labelled node and the total square-average number of loops at a
node.
\end{remark}

\begin{remark} It is worth noting that if we let $\gamma$ evolve from $0$ to $1$, the distribution $F_{\mathcal{AM}_k}(\gamma)$ changes considerably: at $\gamma = 0$, the only path produced with probability $1$ is the one path which has no rises, whereas at $\gamma = 1$, each path is produced with equal probability ($F_{\mathcal{AM}_k}(1)$ is the uniform distribution on alternating Motzkin paths). This phenomenon is reminiscent of percolation processes.
\end{remark}

\begin{example}
\label{example_id3}
For $k = 3$, we compute the average rise-by-altitude vector $\vec{R}$ and the average level-by-altitude vector $\vec{L}$ for alternating Motzkin paths of length 6 as follows. 

\begin{displaymath}
\begin{array}{c@{\hspace{8mm}}c@{\hspace{8mm}}c@{\hspace{8mm}}c@{\hspace{8mm}}c}
\epsfig{file=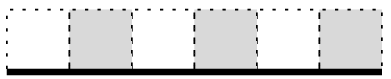, width=0.15\textwidth} & \epsfig{file=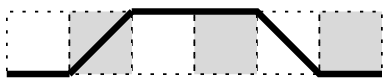, width=0.15\textwidth} & \epsfig{file=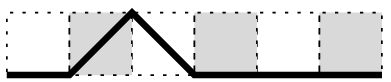, width=0.15\textwidth} & \epsfig{file=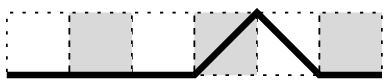, width=0.15\textwidth} & \epsfig{file=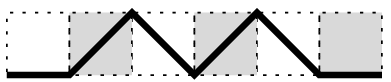, width=0.15\textwidth}\\
\vec{R}=(0,0,0) & \vec{R}=(1,0,0) & \vec{R}=(1,0,0) & \vec{R}=(1,0,0) & \vec{R}=(2,0,0)\\
\vec{L}=(3,0,0) & \vec{L}=(1,1,0) & \vec{L}=(2,0,0) & \vec{L}=(2,0,0) & \vec{L}=(1,0,0)
\end{array}
\end{displaymath}
\begin{displaymath}
E[\vec{R}] = \frac{1}{1+3\gamma+\gamma^2}(3\gamma+2\gamma^2,0,0) \qquad \qquad E[\vec{L}] = \frac{1}{1+3\gamma+\gamma^2}(3+5\gamma+\gamma^2,\gamma,0)
\end{displaymath}

This gives
\begin{eqnarray*}
\| E[\vec{R}]\|^2_2~+~\gamma~\| E[\vec{L}]\|_2^2 & = & \frac{\left ( (3\gamma+2\gamma^2)^2 + \gamma \left((3+5\gamma+\gamma^2)^2+\gamma^2) \right) \right)}{(1+3\gamma+\gamma^2)^2}  \\
& = & \frac{9 \gamma+ 39\gamma^2+44 \gamma^3+14 \gamma^4+\gamma^5}{(3\gamma+2\gamma^2)^2} \\
& = & \frac{N_6(\gamma)}{N_3(\gamma)^2} - 1~.
\end{eqnarray*}
\end{example}

In addition to the three theorems proved in Section \ref{the3}, we give below two more identities involving Catalan and Narayana numbers, for which we do not have combinatorial proofs. These arise as the first-order terms in the asymptotical expansions of the moments of the eigenvalue distribution of $\beta$-Hermite and $\beta$-Laguerre ensembles, and are proved algebraically in \cite{dumitriu03b}. We discuss these in Section \ref{the2}.

\begin{theorem} \label{id4} Using the notations defined above,
\[
\sum_{p \in \mathcal{D}_k} \sum_{i=0}^{k-1} \frac{R_i}{2} (2i+3-R_i) = \sum_{q \in \mathcal{D}_k} \sum_{i=0}^{k-1} {V_i +1 \choose 2}~.   
\]
\end{theorem}

\begin{theorem} \label{id5} Using the notations defined above,
\[
\sum_{p \in \mathcal{AM}_k} \gamma^{r} \left(\sum_{i=0}^{k-1} (i+1) R_i + \gamma \sum_{i=0}^{k-1} i L_i \right) = \sum_{p \in \mathcal{AM}_k} \gamma^{r} \left(\sum_{i=0}^{k-1} {R_i \choose 2} + \gamma \sum_{i=0}^{k-1} {L_i \choose 2} \right)~.
\]
\end{theorem}

\section{The bijection and its variations} \label{the3}

In this section we present one basic construction and three modifications; we use the first two to prove Theorems \ref{id1} and \ref{id2}, and the last two to prove Theorem \ref{id3}.

\subsection{Basic construction} \label{first}

We prove Theorem \ref{id1} by constructing a bijection.

Given an integer $k$, let $p_1$ and $p_2$ be two Dyck paths of length $2k$. Let $i$ be an integer between $0$ and $k-1$, $x_1$ be a rise in $p_1$ from altitude $i$ to altitude $i+1$, and $x_2$ be a fall in $p_2$ from altitude $i+1$ to altitude $i$. To the five-tuplet $(p_1,p_2,i,x_1,x_2)$ we will associate a Dyck path $P$ of length $4k$ which has altitude $2i+2$ in the middle, between steps $2k$ and $2k+1$.

We construct $P$ from $p_1$ and $p_2$ as described below; each move on $p_1$ is followed by a mirror-reversed move in $p_2$, i.e. instead of going left we go right, instead of looking for rises we look for falls and the reverse, instead of flipping up we flip down, etc. 

\begin{figure}[!ht]
\begin{center}
\epsfig{figure = 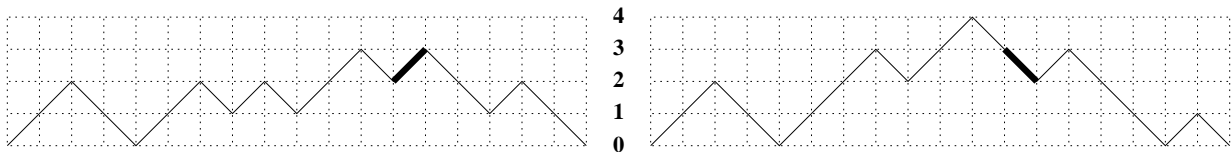, width=\textwidth}
\end{center}
\caption{Choosing a rise $x_1$ from altitude $i=2$ in $p_1$ (left) and a fall $x_2$ from altitude 3 in $p_2$ (right).} \label{BB1}
\end{figure}

\begin{figure}[!ht]
\begin{center}
\epsfig{figure = 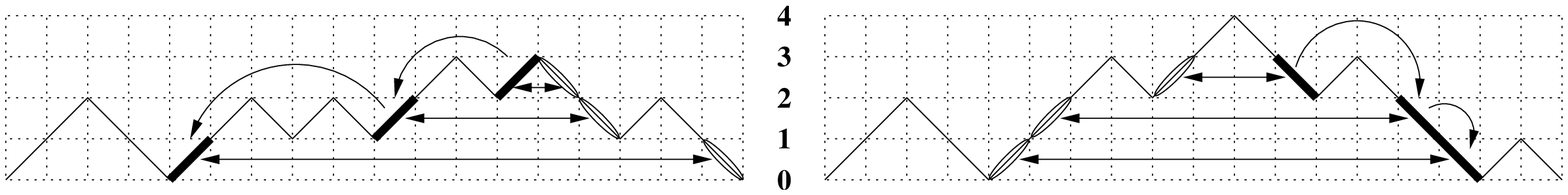, width=\textwidth}
\end{center}
\caption{Finding the ``first rise'' steps from $0$ to $2$ in $p_1$ (left), and the ``first fall'' steps from $2$ to $0$ in $p_2$ (right); the curved arrows point them, and the horizontal double arrows find their respective marked ``closing'' steps.} \label{BB2}
\end{figure}

\fboxsep3pt
\hspace{-2pt}\colorbox[gray]{0.9}{\parbox{\textwidth}{\textit{Step \textbf{1a.}} 
In $p_1$ start at $x_1$, and go left along the path as in Figure \ref{BB1}, the picture on the left, then find the first rise from altitude $i-1$ to altitude $i$, then go left and mark the first rise from $i-2$ to $i-1$, etc. Each of these $i+1$ edges ($x_1$ included) has a ``closing'' fall on the right side of $x_1$, which we find and mark as in the diagram on the left of Figure \ref{BB2}.}}

\hspace{-2pt}\colorbox[gray]{0.9}{\parbox{\textwidth}{\textit{Step \textbf{1b.}} In $p_2$, start at $x_2$, and go right as in the right diagram of Figure \ref{BB1}. Perform the same operations as in Step \textbf{1a}, but mirror-reversed as in the right diagram of Figure \ref{BB2}.}}

\begin{figure}[!ht]
\begin{center}
\epsfig{figure = 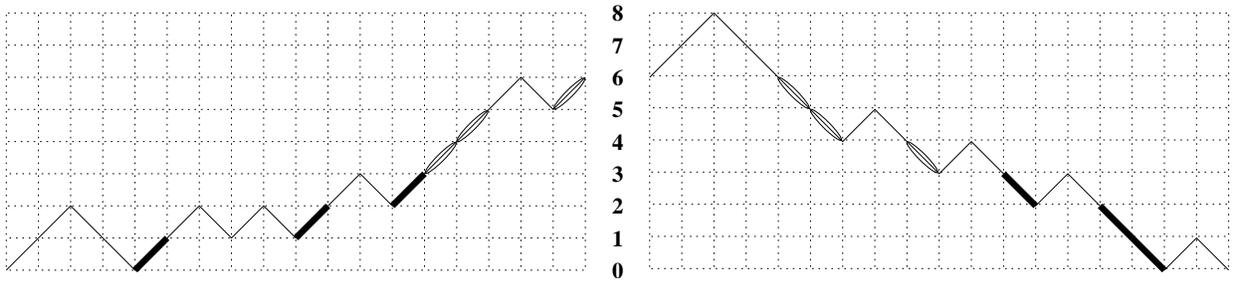, width=\textwidth}
\end{center}
\caption{Flipping the rises in $p_1$ and the falls in $p_2$. The flipped edges correspond to the rightmost rise from altitude $i+1$, the rightmost rise from altitude $i+2$, and so on, in the new path; same for $p_2$ after reversal.} \label{BB3}
\end{figure}

\hspace{-2pt}\colorbox[gray]{0.9}{\parbox{\textwidth}{\textit{Step \textbf{2a.}}
Flip all the closing marked falls in $p_1$ to rises; each flip increases the final altitude of the path by $2$, so the end vertex is at altitude $2i+2$. Note that that the flipped edges correspond, in the new path, to the rightmost rise from altitude $i+1$, the rightmost rise from altitude $i+2$, etc. Hence, given a path of length $2k$ made of $k+i+1$ rises and $k-i-1$ falls which does not go below the $x$-axis, there is a simple transformation which flips the $i+1$ rightmost rises from altitude $i+1$, $i+2$, etc, to falls to get a Dyck path. Thus this process is reversible as demonstrated in Figure \ref{BB3} (on the left).}}

\hspace{-2pt}\colorbox[gray]{0.9}{\parbox{\textwidth}{\textit{Step \textbf{2b.}} Perform the mirror-reversed process on $p_2$, flipping the marked rises to falls; each flip increases the altitude of the initial vertex by $2$, so that at the end, the initial vertex is at altitude $2i+2$.  The process is reversible as demonstrated in Figure \ref{BB3} (on the right).}}

\begin{figure}[!ht]
\begin{center}
\epsfig{figure = 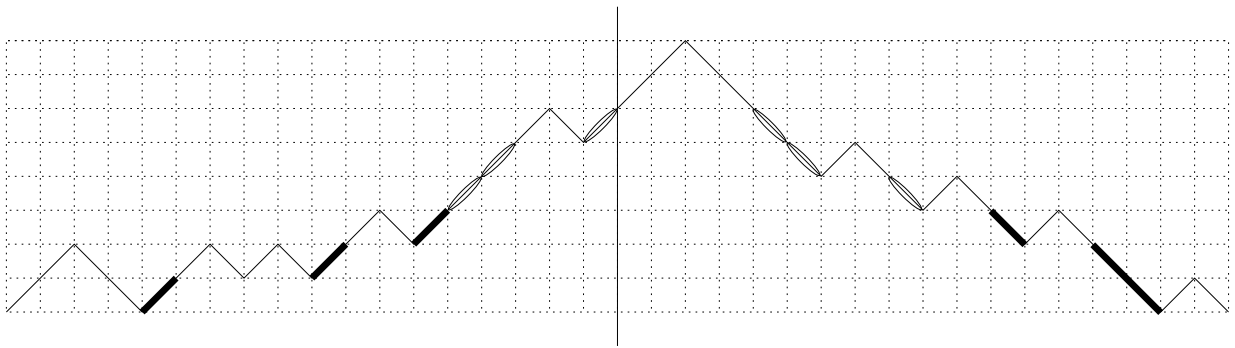, width=\textwidth}
\end{center}
\caption{Concatenating the two paths from Figure \ref{BB3}; the resulting path is a Dyck path of double length and altitude $6=2\times 3$ in the middle.} \label{BB4}
\end{figure}

\hspace{-2pt}\colorbox[gray]{0.9}{\parbox{\textwidth}{\textit{Step \textbf{3.}} We concatenate the two paths obtained from $p_1$ and $p_2$ to obtain a Dyck path of length $4k$ which has altitude $2i+2$ in the middle, between steps $2k$ and $2k+1$, as in Figure \ref{BB4}.}}

The 3-step process above is reversible in a one-to-one and onto fashion. Thus to each five-tuplet $(p_1, p_2, i, x_1, x_2)$ we have associated bijectively a Dyck path $P$ of length $4k$ and altitude $2i+2$ in the middle. 

We can now prove Theorem \ref{id1} merely by counting the two sets described above. 

\textit{Proof of Theorem \ref{id1}.}
Any Dyck path of length $4k$ is at an even altitude in the middle. We separate the Dyck paths which are at altitude $0$ in the middle; since both the left half and the right half of such a path are Dyck paths of length $2k$, it follows that the cardinality of the set 
\[
S_{\textit{right}} = \{P~|~P \in \mathcal{D}_{4k} ~\mbox{and} ~P~\mbox{has positive altitude in the middle} \}
\]
is $|S_{right}| = C_{2k} - C_k^2$.

On the other hand, the cardinality of the set 
\begin{eqnarray*}
S_{\textit{left}}  & = & \{(p_1, p_2, i, x_1, x_2) ~|~p_1 \in \mathcal{D}_{k}, p_2 \in \mathcal{D}_k, i \in \{0,\ldots, k-1\}, \\
& & ~~~~~~~~~~~~~~~~~~~~~~~~~x_1 ~\mbox{a rise at altitude $i$ in $p_1$}, \\
& & ~~~~~~~~~~~~~~~~~~~~~~~~~x_2 ~\mbox{a fall from altitude $i+1$ in $p_2$} \}
\end{eqnarray*} 
is 
\[
S_{\textit{left}} = \sum_{p_1, p_2 \in \mathcal{D}_{k}} \sum_{i=0}^{k-1} R_i(p_1) R_i(p_2)~;
\]
dividing both $S_{\textit{left}}$ and $S_{\textit{right}}$ by $C_{k}^2$ to compute expectations completes the proof. \qed

\subsection{A slight variation} \label{second}

In this section, we slightly modify the construction of Section \ref{first} to make it suitable for the proof of Theorem \ref{id2}.

Given an integer $k$, let $p_1$ and $p_2$ be two Dyck paths of length $2k$. Let $i$ be an integer between $0$ and $k-1$, $x_1$ be a vertex in $p_1$ at altitude $i$, and $x_2$ be a vertex in $p_2$ at altitude $i$. To the five-tuplet $(p_1,p_2,i,x_1,x_2)$ we will associate a Dyck path $P$ of length $4k+2$ which has altitude $2i+1$ in the middle. Note that \emph{all} Dyck paths of length $4k+2$ are at odd altitude in the middle, between steps $2k+1$ and $2k+2$.

Just as before, we construct $P$ from $p_1$ and $p_2$ as described below; each move on $p_1$ is followed by a mirror-reversed move in $p_2$, i.e. instead of going left we go right, instead of looking for rises we look for falls, instead of flipping up we flip down, etc.

We rewrite the construction process below.

\hspace{-2pt}\colorbox[gray]{0.9}{\parbox{\textwidth}{\textit{Step \textbf{1a}.} 
In $p_1$ start at $x_1$, and go left; if $i>0$, find the first rise from altitude $i-1$ to altitude $i$, then go left and mark the first rise from $i-2$ to $i-1$, etc. Each of these $i$ edges has a ``closing'' fall on the right side of $x_1$, which we find and mark. If $i=0$, we mark nothing in the path.}}

\hspace{-2pt}\colorbox[gray]{0.9}{\parbox{\textwidth}{\textit{Step \textbf{1b.}} In $p_2$, start at $x_2$, and go right. Perform the same operations as in Step \textbf{1a}, but mirror-reversed.}}

\hspace{-2pt}\colorbox[gray]{0.9}{\parbox{\textwidth}{\textit{Step \textbf{2a.}}
Flip all the closing marked falls in $p_1$ to rises; each flip increases the final altitude of the path by $2$. \emph{In addition, insert a rise to the right of $x_1$}; the total increase in the altitudeof the end vertex is $2i+1$.}}

Note that  that the inserted edge corresponds in the new path to the rightmost rise from altitude $i$, and the flipped edges correspond to the rightmost rises from altitude $i+1$, $i+2$, etc. Hence, given a path of length $2k+1$ made out of $k+i+1$ rises and $k-i$ falls, which does not go below the $x$-axis, there is a simple transformation which deletes the rightmost rise from altitude $i$ and then flips the $i$ righmost rises from altitude $i$, $i+1$, etc, to falls to get a Dyck path.

\hspace{-2pt}\colorbox[gray]{0.9}{\parbox{\textwidth}{\textit{Step \textbf{2b.}} Perform the mirror-reversed process on $p_2$, flipping the marked rises to falls; each flip increases the initial altitude by $2$. \emph{Add a fall to the left of $x_2$}; the total increase in the altitude of the initial vertex is $2i+1$.}}

\hspace{-2pt}\colorbox[gray]{0.9}{\parbox{\textwidth}{\textit{Step \textbf{3.}} We concatenate the two paths obtained from $p_1$ and $p_2$ to obtain a Dyck path of length $4k+2$ which has altitude $2i+1$ in the middle, between steps $2k+1$ and $2k+2$.}}

The 3-step process above is reversible in a one-to-one and onto fashion. Thus to each five-tuplet $(p_1, p_2, i, x_1, x_2)$ we have associated bijectively a Dyck path $P$ of length $4k+2$ and altitude $2i+1$ in the middle. 

\textit{Proof of Theorem \ref{id2}.}
Once again, we count the sizes of the sets between which we have constructed a bijection; the right set has cardinality $C_{2k+1}$, since any Dyck path of length $4k+2$ has altitude $2i+1$ in the middle, for some $i$. So
\[
S_{\textit{right}} = C_{2k+1}~.
\]

On the other hand, the cardinality of the set 
\begin{eqnarray*}
S_{\textit{left}}  & = & \{(p_1, p_2, i, x_1, x_2) ~|~p_1 \in \mathcal{D}_{k}, p_2 \in \mathcal{D}_k, i \in \{0,\ldots, k\}, \\
& & ~~~~~~~~~~~~~~~~~~~~~~~~~x_1 ~\mbox{a vertex at altitude $i$ in $p_1$}, \\
& & ~~~~~~~~~~~~~~~~~~~~~~~~~x_2 ~\mbox{a vertex at altitude $i$ in $p_2$} \}
\end{eqnarray*} 
is 
\[
S_{\textit{left}} = \sum_{p_1, p_2 \in \mathcal{D}_{k}} \sum_{i=0}^{k} V_i(p_1) V_i(p_2)~.
\]
We then divide both $S_{\textit{left}}$ and $S_{\textit{right}}$ by $C_k^2$ to compute expectations and complete the proof. \qed

\subsection{A version for alternating Motzkin paths} \label{third}

The basic version of the construction and its slight modification presented in Sections \ref{first} and \ref{second} work for Dyck paths; in this section we adapt the construction to work for alternating Motzkin paths. We present two more bijections which we use to prove Theorem \ref{id3}.

To each Motzkin path of length $r$ we will from now on associate a weight $\gamma^r$.

First we need the following lemma.

\begin{lemma}
Let $p$ be an alternating Motzkin path of length $2k$, and $i$ an integer between $0$ and $k-1$. The number of level steps taken in $p$ at altitude $i$ is even, and exactly half of them are on even-numbered steps.
\end{lemma}

\begin{proof}
Let us examine a maximal sequence of level steps at altitude $i$; we use ``maximal'' to express the fact that the steps preceding and succeeding the sequence of level steps (if they exist) are rises or falls. For the benefit of the reader, we include Figure \ref{AltPathParity}.

Assume $i>0$, so that there are steps preceding and succeeding the sequence of level steps.

If the sequence of level steps has even length, then half of them are on even-numbered steps. Moreover, due to the alternating constraint, they have to either be preceded by a rise and succeeded by a fall, or the reverse, as in the regions B, D of Figure \ref{AltPathParity}.

\begin{figure}[ht]
\begin{center} 
\epsfig{figure = 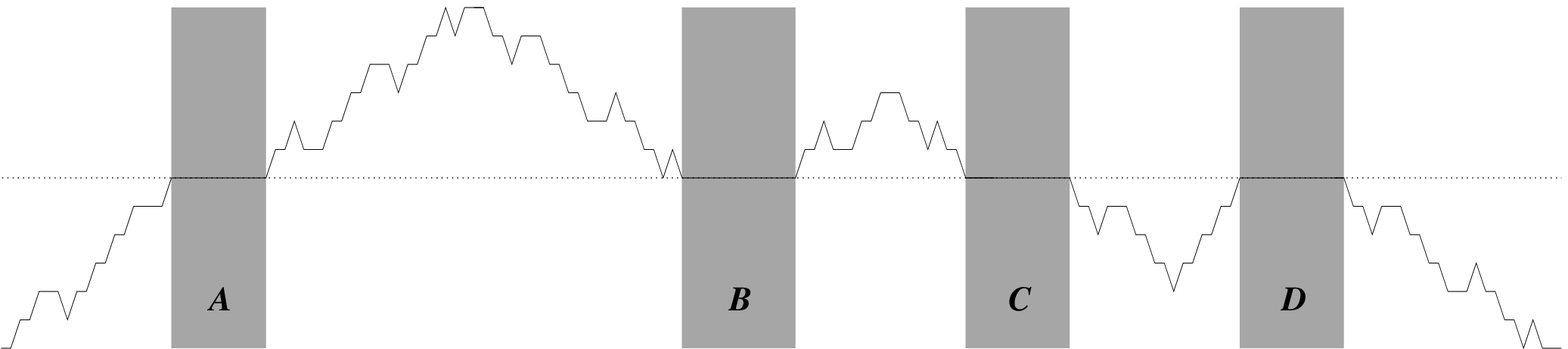, width=\textwidth}
\end{center}
\caption{The four types of maximal sequences of level steps found at some level $i$ in an alternating Motzkin path: even-length ones (B, D) and odd-length ones (A, C).} \label{AltPathParity}
\end{figure}

If the sequence of steps has odd length, there are two possibilities: either both the preceding and the succeeding steps are rises as in region A in Figure \ref{AltPathParity}), or they are both falls as in region C. It is enough to examine the first case (region A). We call the maximal sequence of level steps $a$.

In the first case, the path climbs to a higher altitude, and since it ends at $(2k,0)$, it will have to eventually go below altitude $i$; there is a closest place where the path descends below altitude $i$. The only way in which the path can return to, and then leave, level $i$ is by a sequence \emph{fall, level, level, $\ldots$, level, fall} (see region C). This sequence contains a maximal sequence of level steps (which we call $c$), which has odd length. Moreover, because of the alternating property, the pair of maximal level-$i$ sequences $(a,c)$ will have exactly half of its steps on odd-numbered steps.

Note that the path cannot have two regions of type A without a region of type C between them (nor the converse), since a region of type $A$ implies that a descent to altitude $i-1$ has already taken place and the only way in which this can happen is by passing through a region of type C. So the regions of types A and C alternate in the path, with a region of type A being first and a region of type C being last.

Thus, we can pair \emph{all} the odd-length maximal level sequences at altitude $i$ (each region of type $A$ gets paired with the following region of type C), so that each pair has exactly half of its steps on odd-numbered steps; this shows the claim for $i>0$.

Assume now $i=0$. If there are both preceding and succeeding steps, they can only be a fall and a rise (in this order); in this case the sequence of level steps has even length. Suppose that either the preceding step or the succeeding step is missing (i.e. we are at one end of the path or at the other). 

In the first case the succeeding step can only be a rise, so the path has odd length and one more odd-numbered step than even-numbered steps. We thus know that any alternating Motzkin path starts with an odd-length sequence of level steps. Similarly, it ends with an odd-length sequence of level steps; this sequence has one more even-numbered step. Hence the pair formed by the first and last maximal sequences of level steps at level $0$ has exactly as many odd-numbered steps as even-numbered steps. 

This concludes the proof.
\end{proof}
 
We can now present the new constructions.

Given three integers $k>0$, $k-1\geq r_1, r_2 \geq 0$, let $p_1$ and $p_2$ be two alternating Motzkin paths of length $2k$, with $r_1$ and $r_2$ rises respectively. 

Let $i$ be an integer between $0$ and $k-1$, $x_1$ be a rise in $p_1$ from altitude $i$ to altitude $i+1$, and $x_2$ be a fall in $p_2$ from altitude $i+1$ to altitude $i$. Also, let $y_1$ be a level step in $p_1$ at altitude $i$ which is on an even-numbered step, and $y_2$ a level step in $p_2$ at altitude $i$ which is on an odd-numbered step.

To the five-tuplet $(p_1,p_2,i,x_1,x_2)$ we will associate an alternating Motzkin path $P$ of length $4k$ which has altitude $2i+2$ in the middle, and $r_1+r_2$ rises. 

To the five-tuplet $(p_1, p_2, i, y_1, y_2)$ we will associate an alternating Motzkin path $Q$ of length $4k$ which has altitude $2i+1$ in the middle, and $r_1+r_2+1$ rises.

Just as before, we construct $P$ and $Q$ from $p_1$ and $p_2$ as described below; each move on $p_1$ is followed by a mirror-reversed move in $p_2$.

Note that we no longer can flip rises to falls and vice-versa, since the alternating property would not be respected. 

We present the two constructions below. 

\begin{figure}[ht]
\begin{center}
\epsfig{figure = 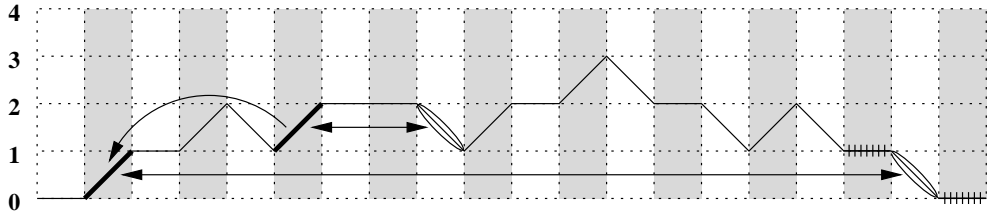, width=0.8\textwidth}
\end{center}
\caption{Choosing the rise; finding the corresponding ascending sequence, the closing one, and the closest level steps. The thick lines represent the ascending sequence, the tripled lines -- the closing one, and the hatched lines are the closest level steps.} \label{AltDyckBijection1}
\end{figure}

\begin{figure}[ht]
\begin{center}
\epsfig{figure = 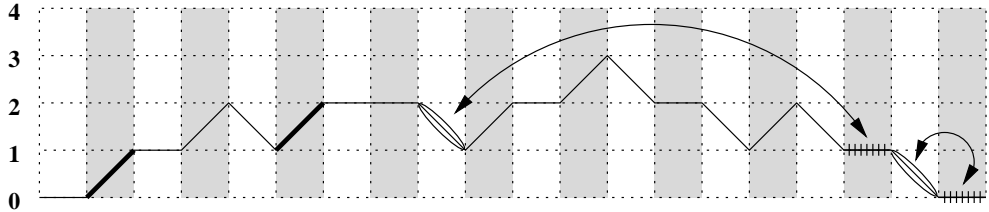, width=0.8\textwidth}
\end{center}
\caption{The marked falls and their corresponding closest level steps to the right.} \label{AltDyckBijection2}
\end{figure}

\hspace{-2pt}\colorbox[gray]{0.9}{\parbox{\textwidth}{\textit{Step \textbf{1a}.} 
In $p_1$ start at $x_1$, and go left; find the first rise from altitude $i-1$ to altitude $i$, then go left and mark the first rise from $i-2$ to $i-1$, etc. Each of these $i+1$ edges has a closing fall on the right side of $x_1$, which we find and mark as in Figure \ref{AltDyckBijection1}. \emph{Each of these marked edges has a closest level step at altitude $i-1$ to the right of it; these are on even-numbered steps. We find them and mark them, as in Figure \ref{AltDyckBijection2}.}}}

\begin{remark} The first descent from $i$ to $i-1$ after the descent from $i+1$ to $i$ must be preceded by a level step, since if it were preceded by a rise, it could not be the first. This is the closest level step at altitude $i$, and it is on an even-numbered step.
\end{remark} 

\hspace{-2pt}\colorbox[gray]{0.9}{\parbox{\textwidth}{\textit{Step \textbf{1b.}} In $p_2$, start at $x_2$, and go right. Perform the same operations as in Step \textbf{1a}, but mirror-reversed.}}

\hspace{-2pt}\colorbox[gray]{0.9}{\parbox{\textwidth}{\textit{Step \textbf{2a.}}
Switch each of the $i+1$ marked falls with its corresponding closest level step to the right, then flip all falls to rises. The new rises are on even-numbered steps, hence the alternance is preserved. Each flip increases the final altitude of the path by $2$, for a total increase in the altitude of the final vertex of $2i+2$, as in Figure \ref{AltDyckBijection2}.}}

Note that to reconstruct the original path, we choose the $i+1$ rightmost rises from $i+1$, $i+2$, etc, find the closest level step to the left of them, switch them and flip the rises to falls to get an alternating Motzkin path.

\hspace{-2pt}\colorbox[gray]{0.9}{\parbox{\textwidth}{\textit{Step \textbf{2b.}} Perform the mirror-reversed process on $p_2$, switching the closest left level steps with the marked rises, and then flipping the marked rises to falls; each flip increases the altitude of the initial vertex by $2$ for a total increase of $2i+2$.}}

\begin{figure}[ht]
\begin{center}
\epsfig{figure = 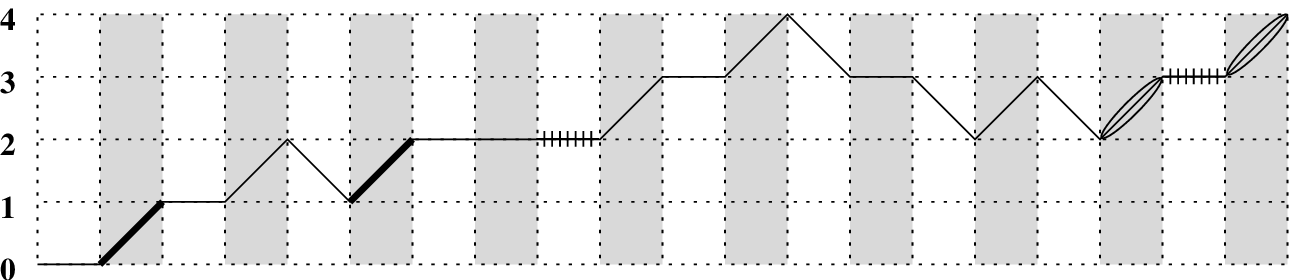, width=0.8\textwidth}
\end{center}
\caption{Switching each closing fall with the corresponding level step, and then flipping the falls to rises.} \label{AltDyckBijection3}
\end{figure}

\hspace{-2pt}\colorbox[gray]{0.9}{\parbox{\textwidth}{\textit{Step \textbf{3.}} We concatenate the two paths obtained from $p_1$ and $p_2$ to obtain an alternating Motzkin path of length $4k$ which has altitude $2i+2$ in the middle.}}

As before, the 3-step process above is reversible in a one-to-one and onto fashion. Thus to each five-tuplet $(p_1, p_2, i, x_1, x_2)$ we have associated bijectively an alternating Motzkin path $P$ of length $4k$ and altitude $2i+2$ in the middle. Moreover, we have not changed the overall number of rises -- we have added $i+1$ rises in $p_1$, but we have flipped $i+1$ rises to falls in $p_2$, and hence the total number of rises in the resulting alternating Motzkin path of length $4k$ is $r_1+r_2$.

The last construction takes a five-tuplet $(p_1, p_2, i, y_1, y_2)$, and produces an alternating Motzkin path of length $4k$ which is at altitude $2i+1$ in the middle. 

The only way in which the last construction differs from the previous one is that it replaces $y_1$ with a rise and $y_2$ with a fall, thus increasing the total number of rises to $r_1+r_2+1$. 
 
\textit{Proof of Theorem \ref{id3}}.
We compute the weight of all the alternating Motzkin paths of length $4k$ which are not at altitude $0$ in the middle; the total weight is 
\[
W_{\textit{left}} = N_{2k}(\gamma) - N_k(\gamma)^2~.
\] 

On the other hand, each alternating Motzkin path which is not at altitude $0$ in the middle is associated to either a five-tuplet $(p_1, p_2, i, x_1, x_2)$, if it is at positive even altitude in the middle, or $(p_1, p_2, i, y_1, y_2)$ if it is at odd altitude in the middle. Hence the total weight can be counted by
\[
W_{\textit{right}} = \sum_{p_1, p_2 \in \mathcal{AM}_k} \gamma^{r_1+r_2} \left(\sum_{i=0}^{k-1} R_i(p_1)R_i(p_2) + \gamma \sum_{i=0}^{k-1} L_i(p_1)L_i(p_2) \right)~,
\]
and dividing by $N_k(\gamma)^2$ to compute expectations, one obtains the statement of Theorem \ref{id3}. \qed

\section{Open problems: two identities} \label{the2}

In this section we present two identities involving Catalan and Narayana numbers which are direct consequences of Theorem 1 in \cite{dumitriu03a}. The proof is algebraic, via random matrix theory. 

Given the nature of the identities, we believe in the existence of a direct combinatorial proof based on constructing a bijection similar to the ones we presented in Section \ref{the3}. We propose the problem of finding such a proof to the interested reader.

\textit{Open Problem 1.} Given a positive integer $k$, let $\mathcal{D}_k$ be the set of all Dyck paths of length $2k$, and given a path $p \in \mathcal{D}_k$, let $\vec{R} = (R_0, \ldots, R_{k-1})$ be the rise-by altitude vector, and $\vec{V} = (V_0, \ldots, V_k)$ be the vertex-by-altitude vector.

Find a combinatorial proof of the following identity (given in Theorem~\ref{id4}):
\[
\sum_{p \in \mathcal{D}_k} \sum_{i=0}^{k-1} ~\frac{R_i}{2} (2i+3-R_i) ~~ =
~~ \sum_{p \in \mathcal{D}_{k-1}} \sum_{i=0}^{k-1} {V_i+1 \choose 2}~~.
\]

For $k=3$, we use Example~\ref{example_id1} and the diagram below.
\begin{displaymath}
\begin{array}{c@{\hspace{30mm}}c}
\epsfig{file=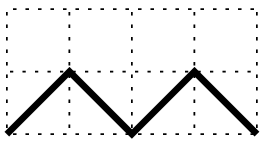} & \epsfig{file=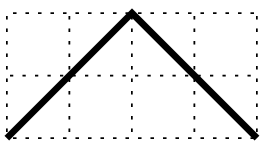}\\[2mm]
\vec{V}=(3,2,0) & \vec{V}=(2,2,1)
\end{array}
\end{displaymath}

From Example~\ref{example_id1}, we deduce that for $k=3$, the sum on the
left is
\begin{eqnarray*}
S_{\textit{left}} & = & \left(\frac{3}{2} ( 3-3)\right) ~~~+ ~~~\left(
\frac{2}{2}(3-2)+\frac{1}{2}(5-1)\right) ~~~+ ~~~\left(
\frac{2}{2}(3-2)+\frac{1}{2}(5-1)\right) ~~~+  \\
& & ~~~+~~~\left(\frac{1}{2}(3-1)+\frac{2}{2}(5-2)\right)
~~~+~~~\left(\frac{1}{2}(3-1)+\frac{1}{2}(5-1)+\frac{1}{2}(7-1)\right)~\\
& = & 16~.
\end{eqnarray*}

On the other hand, from the diagram, the sum on the right is
\[
S_{\textit{right}} ~=~ \left(~{4 \choose 2} ~~+~~{3 \choose 2} ~\right) ~~+~~\left(
~{3 \choose 2} ~~+~~ {3 \choose 2} ~~+~~{2 \choose 2}~\right) ~=~ 16~,
\]
once again.

\textit{Open Problem 2.} Given a positive integer $k$, let $\mathcal{AM}_k$ be the set of all alternating Motzkin paths of length $2k$, and given a path $p \in \mathcal{AM}_k$, let $\vec{R} = (R_0, \ldots, R_{k-1})$ be the rise-by altitude vector, and $\vec{L} = (L_0, \ldots, L_{k-1})$ be the vertex-by-altitude vector.

Find a combinatorial proof of the following identity (given in Theorem~\ref{id5}):
\[
\sum_{p \in \mathcal{AM}_k} \gamma^r \left(\sum_{i=0}^{k-1} (i+1)R_i +
\gamma \sum_{i=0}^{k-1} i L_i \right) ~~=~~ \sum_{p \in \mathcal{AM}_{k-1}}
\gamma^r \left ( \sum_{i=0}^{k-1} {R_i \choose 2} + \gamma \sum_{i=0}^{k-1}
{L_i \choose 2} \right)~~.
\]

We use Example~\ref{example_id3} to illustrate this identity.

The sum on the left is equal to
\begin{eqnarray*}
S_{left} & = & \gamma^0\left(0+0 \gamma\right)
+\gamma^1\left(1+\gamma\right) +\gamma^1\left(1 +0\gamma\right) +
\gamma^1\left(1 +0\gamma\right) + \gamma^2\left(2+0\gamma\right)~,\\
& = & 3(\gamma+\gamma^2)~,
\end{eqnarray*}
while the sum on the right is equal to
\begin{eqnarray*}
S_{right} & = & \gamma^0 \left(0+3\gamma \right) + \gamma^1
\left(0+0\gamma \right) + \gamma^1 \left(0+\gamma \right) +\gamma^1
\left(0+\gamma \right) + \gamma^2(1+0\gamma)~,\\
& = &3(\gamma+\gamma^2)~, \end{eqnarray*}
once again.

\section{Acknowledgements} We would like to thank David Jackson, Richard 
Stanley, Joel Spencer, and especially Alan Edelman for interesting 
conversations about these identities.

\bibliography{bib_ie}
\bibliographystyle{plain}

\end{document}